\documentclass[reqno,11pt]{amsart}
\usepackage{amsmath}
\usepackage{amssymb}
\usepackage{cite}
\usepackage[mathscr]{eucal}
\usepackage{float}
\usepackage[left=34mm,right=34mm,top=35mm,bottom=30mm]{geometry}

\newtheorem{Th}{Theorem}[section]
\newtheorem{Ex}{Example}[section]

\newtheorem{Lemma}{Lemma}[section]
\newtheorem{Prop}{Proposition}[section]
\newtheorem{Cor}{Corollary}[section]

{}

\newcommand{\rb}{\mathbb{R}}
\newcommand{\nb}{\mathbb{N}}

\newcommand{\jc}{\mathcal{J}}
\newcommand{\ic}{\mathcal{I}}
\newcommand{\ac}{\mathcal{A}}

\newcommand{\ck}{\mathfrak{c}}
\newcommand{\fc}{\mathcal{F}}

\newcommand{\pc}{\mathcal{P}}

\newcommand{\zc}{\mathcal{Z}}

\DeclareMathOperator{\Sfg}{S^f_g}

\DeclareMathOperator{\Fin}{Fin}
\DeclareMathOperator{\EU}{EU}
\DeclareMathOperator{\zh}{\mathcal{Z}_h}
\DeclareMathOperator{\zg}{\mathcal{Z}_g}
\DeclareMathOperator{\zgf}{\mathcal{Z}_g(f)}
\DeclareMathOperator{\zagf}{\mathcal{Z}_{ag}(f)}
\DeclareMathOperator{\zfagf}{\mathcal{Z}_{\lfloor ag\rfloor}(f)}

\DeclareMathOperator{\zmghf}{\mathcal{Z}_{\max\{g,h\}}(f)}
\DeclareMathOperator{\zgalf}{\mathcal{Z}_{g_\alpha}(f)}
\DeclareMathOperator{\zhf}{\mathcal{Z}_h(f)}
\DeclareMathOperator{\zf}{\mathcal{Z}(f)}
\DeclareMathOperator{\zlf}{\mathcal{Z}_l(f)}
\DeclareMathOperator{\zl}{\mathcal{Z}_l}
\DeclareMathOperator{\Exh}{Exh}

\newcommand{\bp}{\begin{proof}}
\newcommand{\ep}{\end{proof}}

\begin{document}

\title[Certain observations on weighted density ideals]{CERTAIN OBSERVATIONS ON IDEALS ASSOCIATED WITH WEIGHTED DENSITY USING MODULUS FUNCTIONS}
\author[P. Das  and S. Das]{Pratulananda Das$^*$ and Subhankar Das$^*$\ }
\newcommand{\acr}{\newline\indent}

\address{\llap{*\,}Department of Mathematics, Jadavpur University, Kolkata-700032, West Bengal, India}
\email{pratulananda@yahoo.co.in, subhankarjumh70@gmail.com}


\begin{abstract}
In this article our main object of investigation is the simple modular density ideals $\zgf$ introduced in [Bose et al., Indag. math., 2018] where $g$ is a weight function, more precisely, $g\in G$, $G=\{g:\omega \to [0,\infty):\frac{k}{g(k)}\not\to 0 \text{ and }\:\: g(k)\to \infty \text{ as }\:\:k\to \infty \}$ and $f$ is an unbounded modulus function. We mainly investigate certain properties of these ideals in line of [Kwela et al, J. math. Anal. Appl., 2019]. For an unbounded modulus function $f$ it is shown that there are $1$ or $\ck$ many functions $g\in G$ generating the same ideal $\zgf$. We then obtain certain interactive results involving the sequence of submeasures $\{\phi_k\}_{k\in \omega}$ generating the ideal $\zgf$ and the functions $g,f$. Finally, we present some observations on $\zgf$ ideals related to the notion of increasing-invariance.
\end{abstract}
\maketitle
\smallskip
\noindent{\bf\keywordsname{}:} {Weight function, Simple density ideal, Modulus function, Simple modular density ideal, Erd\"{o}s-Ulam ideal.}

\section{Introduction}
Throughout the article $\omega$ and $\rb$ will stand for the set $\{0,1,\cdots\}$ and the set of real numbers, respectively. For $b,c\in \omega$ with $b\leq c$ by $[b,c]$ we mean the set $\{b,b+1,\dotsc, c\}$ whereas $|C|$ denotes the cardinality of a set $C$. Recall that a non-empty family $\ic\subseteq \pc(\omega)$ is said to be an \textit{ideal} on $\omega$ if $(1)$ for $B, C\in \ic$, $B\cup C\in \ic$, $(2)$ for $B\in \ic$ and $C\subseteq B$, $C\in \ic$, $(3)$ $\omega\not\in \ic$. As in practice $\Fin$ denotes the ideal on $\omega$ consisting of all finite subsets of $\omega$. An ideal $\ic$ is called free (or admissible) if $\Fin \subset \ic$. Throughout we assume that $\ic$ is an admissible ideal on $\omega$.

Let us now recall the notions of natural density function and density zero ideal of $\omega$ which have been studied extensively over the years having immense importance both in number theory and in set theory. For $C\subseteq \omega$, $\overline{d}(C)$ and $\underline{d}(C)$ denote the \textit{upper} and \textit{lower} natural densities of $C$ (\cite{buck-1, buck-2}), respectively. For $k\in \omega$
\[\overline{d}(C)=\limsup_{k\to \infty}\frac{|C\cap [0,k-1]}{k}\:\: \text{and}\:\: \underline{d}(C)=\liminf_{k\to \infty}\frac{|C\cap [0,k-1]}{k}.\]\\
When $\overline{d}(C)=\underline{d}(C)$, it is said that the \textit{natural density} $d(C)$ exists and $d(C)=\overline{d}(C)=\underline{d}(C)$. The family $\zc=\{C\subseteq \omega:\overline{d}(C)=0\}$ is known as the \textit{classical density zero} ideal. This ideal has been extensively investigated over the years in set theory (one can see the seminal papers \cite{hh,f} for more information).

In the recent past the classical density zero ideal has been generalized in different ways. In \cite{sbpdsp} a generalized form $d_\alpha$ of natural density was introduced. For $C\subseteq \omega$ and $0<\alpha\leq 1$, one defines  $\overline{d}_\alpha(C)=\limsup_{k\to \infty}\frac{|C\cap [0,k-1]}{k^\alpha}$ while $\underline{d}_\alpha(C)$ is defined similarly. The notion was further extended to a more generalized form $\overline{d}_g$ \cite{mbpdmfjs} using a weight function $g$, where $g\in G$ and $G$ stands for all functions $g: \omega \to [0, \infty)$ satisfying the property that $g(k) \to \infty$  and $\frac{k}{g(k)}\not\to 0$ as $k\to \infty$. To be precise, for $C\subseteq \omega$, $$\overline{d}_g(C)=\limsup_{k\to \infty}\frac{|C\cap [0,k-1]|}{g(k)}.$$ This density function has since then come to be known as \textit{simple density} function and the associated ideal  $\zg=\{C\subseteq \omega:\overline{d}_g(C)=0\}$ as the \textit{simple density} ideal (the term coined in \cite{ak-2}). Subsequently, different aspects of simple density ideals have been deeply investigated in \cite{ak-1, ak-2}.

We now come to the main focal theme of this article. In \cite{kbpdak} using a modulus function a generalized form of the simple density ideal was introduced  which was called the \textit{modular simple density} ideal $\zgf$, where $g\in G$ and $f$ is a modulus function (In \cite{jt}, this ideal was called \textit{simple generalised density ideal}). This ideal has subsequently been investigated to some extent in \cite{pdag}. It has been observed that the ideal $\zgf$ generally inherits similar kinds of properties which $\zg$ possess. Continuing in this line, in this article we further explore various properties of $\zgf$ which shows that many of the interesting results of the articles \cite{ak-1, ak-2} can be extended to this more general settings.

The article is arranged as follows. In Section 2, definitions of some well known notions and some basic observations on $\zgf$ ideals are presented. In particular it is observed that $\cap_{g\in G}\zgf=\Fin$ and $\cup_{g\in G}\zgf=\zl$ for any modulus function $f$. In Section 3, we are able to show that for an unbounded modulus function there are $1$ or $\ck$ many functions $g$ generating the same ideal $\zgf$. In Section 4 we present certain interactive results related to the sequence of submeasures $\{\phi_k\}_{k\in \omega}$ generating the ideal $\zgf$ and the functions $g,f$. It is shown that $\sup_{k\in \omega} \phi_k(\omega)<\infty$ is equivalent to the fact that the sequence $\{\frac{f(k)}{f(g(k))}\}_{k\in \omega}$ is bounded. Finally, some observations on $\zgf$ ideals related to the notion of increasing-invariance are presented.

\section{Preliminaries}

A \textit{submeasure} on $\omega$ is a function $\phi:\pc(\omega)\to [0,\infty)$, where $\phi$ satisfies the conditions $(1)$ $\phi(\emptyset)=0$, $(2)$ for $B,C\subseteq \omega$ and $B\subseteq C$, $\phi(B)\leq \phi(C)$, $(3)$ $\phi(B\cup C)\leq \phi(B)+\phi(C)$, $(4)$ for all $n\in \omega$ $\phi(\{n\})<\infty$. When $\phi(C)=\lim_{k\to \infty}\phi(C\cap [0,k-1])$ for all $C\subseteq \omega$, $\phi$ is said to be a \textit{lower semicontinuous submeasure} (in short, lscsm). For a lscsm $\phi$, $\Exh(\phi)$ is an exhaustive ideal on $\omega$ and $\Exh(\phi)=\{C\subseteq \omega:\lim_{k\to \infty}\phi(C\setminus [0,k-1])=0\}$. It is well known that $\zc=\Exh(\phi)$, where $\phi(C)=\sup_{k\in \omega}\frac{|C\cap [0,k-1]|}{k}$ for $C\subseteq \omega$ (for details one should consult the seminal papers by Solecki \cite{solecki, solecki1}).

Coming now to the definition of simple density, for $g\in G$ $\zg\neq \{\emptyset\}$ is equivalent to $g(k)\to \infty$ as $k\to \infty$ and $\omega\not\in \zg$ is equivalent to $\frac{k}{g(k)}\not\to 0$ as $k\to \infty$ (see \cite{mbpdmfjs}).  $\zc_{id_\omega}=\zc$, where $id_{\omega}$ is the identity function on $\omega$. Also, $\zg=\zc_{\lfloor g \rfloor}$ for $g\in G$. Moreover, $H=\omega^\omega\cap G$ and $H^\uparrow=\{g\in H:g\:\: \text{is nondecreasing}\}$ (see \cite{ak-1}).

A function $f:\rb^+\cup \{0\}\to \rb^+\cup \{0\}$ is called a \textit{modulus} function (see \cite{agb}) if the following conditions hold $(1)$ $f(x)=0$ if and only if $x=0$, $(2)$ for all $x,y\in \rb^+\cup \{0\}$ $f(x+y)\leq f(x)+f(y)$, $(3)$ $f$ is right continuous at $0$, $(4)$ $f$ is increasing. Some  examples of modulus functions are $(1)$ $f(x)=x$, $x\in \rb^+\cup \{0\}$, $(2)$ $f(x)=\frac{x}{1+x}$, $x\in \rb^+\cup \{0\}$, $(3)$ $f(x)=x^\beta$, $x\in \rb^+\cup \{0\}$, for $\beta\in (0,1)$, $(4)$ $f(x)=log(1+x)$, $x\in \rb^+\cup \{0\}$.
$\overline{d}^f(C)$ denotes the \textit{upper $f$ density} function (see \cite{agb}), where $\overline{d}^f(C)=\limsup_{k\to \infty}\frac{f(|C\cap [0,k-1]|)}{f(k)}$ for $C\subseteq \omega$. $\underline{d}^f$ (the \textit{lower $f$ density} function) is similarly defined. It is easy to note that $\zf=\{C\subseteq \omega:\overline{d}^f(C)=0\}$ forms ideal on $\omega$ (see \cite{agb}) (though this ideal has not been studied on its own, but one should keep in mind that several of its properties come as special cases of the properties of the ideals $\zgf$ taking $g$ as the identity function). By $\zlf$, we mean the family $\{C\subseteq \omega:\underline{d}^f(C)=0\}$. Another set naturally associated with the natural density function is the set $\zl=\{C\subseteq \omega:\underline{d}(C)=0\}$. It is quite clear that $\zlf\subseteq \zl$.

Let $f$ be an unbounded modulus function.  One can naturally think of uniting the two approaches mentioned earlier, one by using a weight function and the other using a modulus function by defining $$\overline{d}_g^f(C)=\limsup_{k\to \infty}\frac{f(|C\cap [0,k-1]|)}{f(g(k))}~~ \mbox{and}~~  \underline{d}_g^f(C)=\liminf_{k\to \infty}\frac{f(|C\cap [0,k-1]|)}{f(g(k))}.$$ If $\overline{d}_g^f(C)=\underline{d}_g^f(C)$, then $d_g^f(C)$ exists and $d_g^f(C)=\overline{d}_g^f(C)=\underline{d}_g^f(C)$. $d_g^f$ is called the \textit{modular simple density} function (see \cite{kbpdak}). The corresponding \textit{modular simple density} ideal is $\zgf=\{C\subseteq \omega:d_g^f(C)=0\}$ (see \cite{kbpdak}). Note that $\zgf=\zg$ when $f$ is the identity function.

Next we recall certain notions from the theory of ideals which will be useful in the sequel. Let $\{I_k\}$ be a sequence of pairwise disjoint and finite intervals and let $\{\phi_k\}$ be a sequence of submeasures concentrated on $I_k$. The ideal $\ic=\{C\subseteq \omega:\lim_{k\to \infty}\sup_{k\in \omega}\phi_k(C)=0\}$ is called a generalised density ideal (see \cite{f-1}). Let $f:\omega\to [0,\infty)$ be a function for which $\Sigma^\infty_{i=0}f(i)=\infty$ and $\lim_{k\to \infty}\frac{f(k)}{\Sigma^k_{i=0}f(i)}=0$. Then $\ic$ is an Erd\"{o}s-Ulam ideal (in short, $\EU$ ideal) if $\ic=\{C\subseteq \omega:\lim_{k\to \infty}\frac{\Sigma_{i\in [0,k-1]\cap C}f(i)}{\Sigma^{k-1}_{i=0}f(i)}=0\}$ (see \cite{just}).
$\ic$ is \textit{increasing-invariant} if for each $C\in \ic$ and $B\subseteq \omega$ satisfying $|B\cap [0, k-1]|\leq |C\cap [0,k-1]|$ for all $k\in \omega$, we have $B\in \ic$. If $\ic$ is increasing invariant and $f:\omega\to \omega$ is injective, increasing and $C\in \ic$, then $f(C)\in \ic$ (see \cite{mbsgjs, ak-2}). If $C\not\in \ic$, then $\ic|C=\{B\cap C:B\in \ic\}$ is an ideal on $C$. $\ic|C$ is called the restriction of $\ic$ to $C$. Two ideals $\ic$ and $\jc$ on infinite countable sets $X$ and $Y$, respectively, are isomorphic ($\ic\cong \jc$) if there exists a bijection $\phi:X\to Y$ for which $C\in \ic$ if and only if $\phi(C)\in \jc$ for $C\subseteq X$. $\ic\sqsubseteq \jc$ if there exists a bijection $\phi: X\to Y$ for which $C\in \ic$ implies that $\phi(C)\in \jc$ for each $C\subseteq X$. Let $\ic$ and $\jc$ be two ideals on $X$. Then $\ic\leq_K \jc$ (\textit{Kat\v{e}tov order}) if there exists a function $\phi: X\to X$ for which $C\in \ic$ implies that $\phi^{-1}(C)\in \jc$ for each $C\subseteq \omega$. The family $H(\ic)=\{C\subseteq X:\ic|C\cong \ic\}$ is the \textit{homogeneity} family of $\ic$ (see \cite{akjt}). The symbol $id_\omega$ denotes the identity function from $\omega$ onto $\omega$. \\

We end this section with some basic observations on $\zgf$. From \cite[Proposition 2]{ak-1} it is known that $\cap_{g\in G}\zg=\Fin$ and $\cup_{g\in G}\zg=\zl$. Analogously, we can prove that $\cap_{g\in G}\zgf=\Fin$ and $\cup_{g\in G}\zgf=\zlf$ for an unbounded modulus function $f$.

\begin{Prop}
\label{P1}
Let $g\in G$ and $f$ be an unbounded modulous function. The following statements hold.\\
$(1)$ $\cap_{g\in G}\zgf=\Fin$.\\
$(2)$ $\cup_{g\in G}\zgf=\zlf$.
\end{Prop}
\bp
$(1)$. The result easily follows from the fact that $\zgf\subseteq \zg$ for any $g\in G$ (\cite[Proposition 2.6]{kbpdak}) and $\cap_{g\in G}\zg=\Fin$.

$(2)$. Let $C\in \zlf$. Clearly, $\liminf_{k\rightarrow\infty}\frac{f(|C\cap[0,k-1]|)}{f(k)}=0$. Choose a sequence $\{k_m\}_{m\in \omega}$ for which $\lim_{m\rightarrow\infty}\frac{f(|C\cap[0,k_m-1]|)}{f(k_m)}=0$. Consider a $g:\omega\rightarrow [0,\omega)$ defined by $g(k)=\min\{k_m:k\leq k_m\}$. Now $g(k)\rightarrow \infty$ as $k\rightarrow \infty$ and also $\frac{k}{g(k)}\nrightarrow 0$ as $\frac{k_m}{g(k_m)}=1$ for all $m\in \omega$. Take $k\in \omega$ such that $g(k)=k_m$ for some $m\in \omega$ and $k\leq k_m$. Note that

\[\frac{f(|C\cap [0,k-1]|)}{f(g(k))}=\frac{f(|C\cap [0,k-1]|)}{f(k_m)}\leq \frac{f(|C\cap [0,k_m-1]|)}{f(k_m)}.\]

Since $\lim_{m\rightarrow\infty}\frac{f(|C\cap[0,k_m-1]|)}{f(k_m)}=0$ and $k\rightarrow \infty$ implies that $m\rightarrow \infty$, we must have $\lim_{k\rightarrow \infty}\frac{f(|C\cap [0,k-1]|)}{f(g(k))}=0$. Therefore, $C\in \zgf$. Hence, $\zlf\subseteq \cup_{g\in G}\zgf$.

Conversely, let $C\in \cup_{g\in G}\zgf$ and let $C\in \zgf$ for some $g\in G$. As $\omega\not\in \zgf$, $\limsup_{k\rightarrow \infty}\frac{f(k)}{f(g(k))}>0$. Choose an increasing sequence $\{k_m\}_{m\in \omega}$ for which the sequence $\{\frac{f(g(k_m))}{f(k_m)}\}_{m\in \omega}$ is bounded. Also, choose $M>0$ satisfying $\frac{f(g(k_m))}{f(k_m)}\leq M$ for each $m\in \omega$. Now
\[\frac{f(|C\cap[0,k_m-1]|)}{f(k_m)}=\frac{f(|C\cap[0,k_m-1]|)}{f(g(k_m))}\frac{f(g(k_m))}{f(k_m)}\leq \frac{f(|C\cap[0,k_m-1]|)}{f(g(k_m))}M.\]
As $C\in \zgf$, $\lim_{m\rightarrow \infty}\frac{f(|C\cap[0,k_m-1]|)}{f(g(k_m))}=0$. Consequently, $\lim_{m\rightarrow \infty}\frac{f(|C\cap[0,k_m-1]|)}{f(k_m)}=0$. Thus $C\in \zlf$ and so $\cup_{g\in G}\zgf\subseteq \zlf$. Hence, $\cup_{g\in G}\zgf=\zlf$.
\ep

As $\zlf\subseteq \zl$, we have $\cup_{g\in G}\zgf\subseteq \zl$ for any modulus function $f$. In the following example we prove that $\cup_{g\in G}\zgf\varsubsetneqq \zl$ for certain modulus function $f$.
\begin{Ex}
\label{E}
Consider the modulus function $f(x)=log(1+x)$. We prove that $\zlf\varsubsetneqq \zl$. For $m\in \omega$ $\lfloor \sqrt{m}\rfloor\leq \sqrt{m}<m$. Now construct a set $C\subseteq \omega$ which contains exactly $\lfloor \sqrt{m}\rfloor$ elements from $\{0,1,\cdots,m-1\}$ for each $m\in \omega$. We now have
\[\frac{|C\cap [0,m-1]|}{m}= \frac{\lfloor \sqrt{m}\rfloor}{m}\leq \frac{\sqrt{m}}{m}=\frac{1}{\sqrt{m}}.\]
Therefore, $\liminf_{m\to \infty}\frac{|C\cap [0,m-1]|}{m}=0$ and $C\in \zl$. Again,
\[\frac{f(|C\cap [0,m-1]|)}{f(m)}=\frac{log(1+\lfloor \sqrt{m}\rfloor)}{log(1+m)}>\frac{log(1+\sqrt{m}-1)}{log(1+m)}=\frac{\frac{1}{2}log(m)}{log(1+m)}\]
Thus $\lim_{m\to \infty}\frac{f(|C\cap [0,m-1]|)}{f(m)}>\frac{1}{2}$ and $C\not\in \zlf$. Hence, $\zlf\varsubsetneqq \zl$. Consequently, $\cup_{g\in G}\zgf\varsubsetneqq \zl$.
\end{Ex}

Following \cite[Lemma 5]{ak-1}, we have the next result.
\begin{Lemma}
\label{LO1}
Let $g\in G$ and $f$ be an unbounded modulus function. The following statements are equivalent.\\
$(1)$ $\zgf\subseteq \zf$.\\
$(2)$ $\liminf_{k\to \infty}\frac{f(k)}{f(g(k))}>0$.
\end{Lemma}
\bp
$(1)\Rightarrow (2)$. Assume that $\liminf_{k\to \infty}\frac{f(k)}{f(g(k))}=0$. Choose an increasing sequence $\{k_m\}_{m\in \omega}$ for which $\lim_{m\to \infty}\frac{f(k_m)}{f(g(k_m))}=0$ and $k_{m+1}>2k_m$. Let $C=\cup_{m\in \omega}[k_m,2k_m)$. Now $\frac{f(|C\cap [0,2k_m-1]|)}{f(2k_m)}\geq \frac{f(k_m)}{2f(k_m)}=\frac{1}{2}$ (as $f(2k_m)\leq 2f(k_m)$). Therefore, $C\not\in \zf$. For $k\in \omega$ choose $m\in \omega$ with $k_m\leq k<k_{m+1}$. Note that

\[\frac{f(|C\cap [0,k-1]|)}{f(g(k))}\leq \frac{f(|C\cap [0,k-1]|)}{f(g(k_m))}<\frac{f(k_m+|C\cap [k_m,2k_m)|)}{f(g(k_m))}\] \[\leq \frac{f(k_m)+f(|C\cap [k_m,2n_m)|)}{f(g(k_m))}\leq \frac{2f(k_m)}{f(g(k_m))}.\]
Now $k\to \infty$ implies that $m\to \infty$ and $\lim_{m\to \infty}\frac{f(k_m)}{f(g(k_m))}=0$. Thus $C\in \zgf$ and by $(1)$, $C\in \zf$ which is a contradiction. Hence, $\liminf_{k\to \infty}\frac{f(k)}{f(g(k))}>0$.

$(2)\Rightarrow (1)$. Let $C\in \zgf$ and $\varepsilon>0$ be given. Let $\alpha=\liminf_{k\to \infty}\frac{f(k)}{f(g(k))}$. Choose $k_0,k_1\in \omega$ for which $\frac{f(k)}{f(g(k))}>\frac{\alpha}{2}$ for all $k\geq k_0$ and $\frac{f(|C\cap [0,k-1]|)}{f(g(k))}<\frac{\varepsilon\alpha}{2}$ for all $k\geq k_1$. Let $k_2=\max\{k_0,k_1\}$. For $k\geq k_2$ we then obtain
\[\frac{f(|C\cap [0,k-1]|)}{f(k)}=\frac{1}{\alpha}\frac{f(|C\cap [0,k-1]|)}{\frac{f(k)}{\alpha}}<\frac{1}{\alpha}\frac{f(|C\cap [0,k-1]|)}{\frac{f(g(k))}{2}}<\frac{2}{\alpha}\frac{\varepsilon\alpha}{2}=\varepsilon.\]

Therefore, $C\in \zf$. Hence, $\zgf\subseteq \zf$.
\ep


\section{Generating the same $\zgf$ ideal with different weight functions $g$}
Suppose that $f$ is an unbounded modulus function. We know that for $c_1, c_2>0$ and $g,h\in G$ if $c_1\leq \frac{f(g(n))}{f(h(n))}\leq c_2$, then $\zgf=\zhf$ \cite[Proposition 2.3]{kbpdak}. From the above result it clearly follows that for a given $g\in G$ $\zgf=\zagf$ for $a\geq 1$. This implies that  there are $\ck$ many functions which generate the ideal $\zgf$ for a given $g\in G$. Clearly, for each $a\in [1,2]$ $\zgf=\zfagf$. For $g\in H^{\uparrow}$ let $\Sfg=\{h\in H^\uparrow: \zgf=\zhf\}/R$, where $R$ is an equivalence relation defined by $gRh\Leftrightarrow \{\limsup_{n \to \infty}\frac{f(g(n))}{f(h(n))}<\infty\wedge \limsup_{n \to \infty}\frac{f(h(n))}{f(g(n))}<\infty\}$. We now look into the cardinality of $\Sfg$ for a given $g\in H^{\uparrow}$ and a given unbounded modulus function $f$. First we present the following results.

\begin{Lemma}
\label{LS1}
Let $g\in H^{\uparrow}$ and $f$ be an unbounded modulus function. Let $C\subseteq \omega$ and $C=\{c_0<c_1<\cdots\}$. The following statements are equivalent.\\
$(1)$ $C\in \zgf$.\\
$(2)$ $\lim_{m\to \infty}\frac{f(|C\cap [0,c_m-1]|)}{f(g(c_m))}=0$.
\end{Lemma}
\bp
The implication $(1)\Rightarrow (2)$ easily follows from the definition of $\zgf$.

$(2)\Rightarrow (1)$. Let $k\in \omega$ and $k\geq c_0$. Clearly, $k\in [c_m,c_{m+1})$ for exactly one $m\in \omega$. We now have
\[\frac{f(|C\cap [0,k-1]|)}{f(g(k))}\leq \frac{f(|C\cap [0,c_{m+1}-1]|)}{f(g(c_m))}=\frac{f(|C\cap [0,c_m-1]|+1)}{f(g(c_m))}\] \[\leq  \frac{f(|C\cap [0,c_m-1]|)}{f(g(c_m))}+\frac{f(1)}{f(g(c_m))}.\]

Since $f$ is unbounded, $g$ is nondecreasing and by $(2)$, $\lim_{k\to \infty}\frac{f(|C\cap [0,k-1]|)}{f(g(k))}=0$. Hence, $C\in \zgf$.
\ep

\begin{Lemma}
\label{LS2}
Let $g,h\in H^{\uparrow}$ and $f$ be an unbounded modulus function. Suppose that $\zgf=\zhf$. Then $\zgf=\zmghf$.
\end{Lemma}
\bp
Suppose that $C\in \zgf$. For $k\in \omega$ observe that $f(\max\{g(k),h(k)\})=\max\{f(g(k)),f(h(k))\}$ as $f$ is increasing. Then
\[\frac{f(|C\cap [0,k-1]|)}{f(\max\{g(k),h(k)\})}= \frac{f(|C\cap [0,k-1]|)}{\max\{f(g(k)),f(h(k))\}}\leq \frac{f(|C\cap [0,k-1]|)}{f(g(k))}.\]
Consequently, $\lim_{k\to \infty} \frac{f(|C\cap [0,k-1]|)}{f(\max\{g(k),h(k)\})}=0$ and $C\in \zmghf$. Therefore, $\zgf\subseteq \zc_{\max\{g,h\}}(f)$.

Let $C\in \zc_{\max\{g,h\}}(f)$. Put $C_1=\{k\in C:h(k)\leq g(k)\}=\{c_0<c_1<\cdots\}$. Clearly, $C_1\in \zc_{\max\{g,h\}}(f)$. Let $C_2=C\setminus C_1$. Also, $C_2\in \zc_{\max\{g,h\}}(f)$. Now \[\frac{f(|C_1\cap [0,c_m-1]|)}{f(g(c_m)}=\frac{f(|C_1\cap [0,c_m-1]|)}{\max\{f(g(c_m),f(h(c_m))\}}.\] As $C_1\in \zc_{\max\{g,h\}}(f)$, we obtain that $\lim_{m\to \infty}\frac{f(|C_1\cap [0,c_m-1]|)}{f(g(c_m))}=0$. By Lemma \ref{LS1}, $C_1\in \zgf$. Similarly, we have $C_2\in \zhf$. Since $\zgf=\zhf$, $C_2\in \zgf$. Consequently, $C=C_1\cup C_2\in \zgf$ and so $\zmghf\subseteq \zgf$. Hence, $\zgf=\zmghf$.
\ep

For the next two results we follow the lines of argument of \cite[Proposition 5]{ak-1} and \cite[Theorem 1]{ak-1}, respectively.
\begin{Prop}
\label{PS1}
Let $g\in H^\uparrow$ and $f$ be an unbounded modulus function. If $|\Sfg|>1$, then $|\Sfg|=\ck$.
\end{Prop}
\bp
Let $h\in \Sfg$ and let $h\neq g$. Assume that $\limsup_{k\to \infty}\frac{f(h(k))}{f(g(k))}=\infty$. Choose an increasing sequence $\{k_m\}_{m\in \omega}$ for which $\lim_{m\to \infty}\frac{f(h(k_m))}{f(g(k_m))}=\infty$. By an induction process, we may assume that $g(k_{m+1})\geq h(k_m)$ for all $m\in \omega$. Let $\alpha=(\alpha_1,\alpha_2,\cdots)\in 2^\omega$. We now define $g_\alpha:\omega\rightarrow [0,\infty)$ by
\[g_\alpha(k)=\\
\begin{cases}
g(k)& \:\: \text{if}\:\: k< k_0\\

g(k)& \:\: \text{if}\:\: k\in [k_m,k_{m+1}), \alpha_m=0\\

\max\{g(k),h(k)\}& \:\: \text{if}\:\: k\in [k_m,k_{m+1}),\alpha_m=1

\end{cases}\]

Clearly, $g_\alpha\in H^\uparrow$. We show that $\zgalf=\zgf$. Let $C\in \zgf$. For $k\in \omega$ $g(k)\leq g_\alpha(k)\leq \max\{g(k),h(k)\}$ and $f$ is increasing, so $\zgf\subseteq \zgalf\subseteq \zmghf$.
Now $h\in \Sfg$ and so $\zgf=\zhf$. By Lemma \ref{LS2}, $\zgf=\zmghf$. Therefore, from $\zgf\subseteq \zgalf\subseteq \zmghf$ we have $\zgalf=\zgf$. Choose $\alpha=(\alpha_1,\alpha_2,\cdots)$ and $\beta=(\beta_1,\beta_2,\cdots)\in 2^\omega$. When for infinitely many $k\in \omega$ $\alpha_k>\beta_k$, we have $\limsup_{k\to \infty}\frac{f(g_\alpha(k))}{f(g_\beta(k))}=\infty$ and when for infinitely many $k\in \omega$ $\alpha_k<\beta_k$, we have $\limsup_{k\to \infty}\frac{f(g_\beta(k))}{f(g_\alpha(k))}=\infty$. Now there exists a family $\ac\subseteq 2^\omega$ for which $|\ac|=\ck$ and $|\{k\in \omega:\alpha_k=\beta_k\}|$ is finite for all $\alpha,\beta\in \ac$. Hence, $|\Sfg|=\ck$.

\ep

\begin{Th}
\label{TS1}
Let $g\in H^{\uparrow}$ and $f$ be an unbounded modulus function. If $|\Sfg|=1$, then there exist $M>0$ and $\varepsilon>0$ for which $\frac{f(g(k+\lfloor \varepsilon f(g(k)) \rfloor))}{f(g(k))}\leq M$ for all but finitely many $k\in \omega$.
\end{Th}
\bp
Assume the contrary. Then for every $M>0$ and $\varepsilon>0$, $\frac{f(g(k+\lfloor \varepsilon f(g(k)) \rfloor))}{f(g(k))}> M$ for infinitely many $k\in \omega$. We claim that $|\Sfg|>1$.
Choose a sequence $\{k_m\}_{m\in \omega}$ for which $\frac{f(g(k_m+\lfloor \frac{f(g(k_m))}{2^m}\rfloor))}{f(g(k_m))}>m$ and $k_{m+1}>k_m+\frac{f(g(k_m))}{2^m}$. Clearly, $\lfloor \frac{f(g(k_m))}{2^m} \rfloor\leq \frac{f(g(k_m))}{2^m}<k_{m+1}-k_m$ and $\frac{f(g(k_{m+1}))}{f(g(k_m))}=\frac{f(g(k_m+k_{m+1}-k_m))}{f(g(k_m))}>\frac{f(g(k_m+\lfloor \frac{f(g(k_m))}{2^m}\rfloor))}{f(g(k_m))}>m$, $m\geq 2$. Let $I_m=[k_m,k_m+\lfloor \frac{f(g(k_m))}{2^m}\rfloor]$ for each $m$. Choose a function $h$ as follows.
\[h(k)=\\
\begin{cases}
g(k_m+\lfloor \frac{f(g(k_m))}{2^m}\rfloor)& \:\: \text{if}\:\: k\in I_m \:\: \text{for some}\:\: m\in \omega\\

g(k)& \:\: \text{otherwise.}

\end{cases}
\]

Clearly, $h\in H^\uparrow$. As $\frac{f(h(k_m))}{f(g(k_m))}>m$ for all $m\in \omega$, $\limsup_{k\to \infty}\frac{f(h(k))}{f(g(k))}=\infty$. Again, from $h(k)\geq g(k)$ for each $k\in \omega$, it follows that $\zgf\subseteq \zhf$.

We now prove that $\zhf\subseteq \zgf$. Let $C\in \zhf$. Clearly, $C\setminus \cup_{m\in \omega}I_m\in \zhf$ and $C\setminus \cup_{m\in \omega}I_m\subseteq \{k\in \omega:h(k)=g(k)\}$ enumerate the elements of $C\setminus \cup_{m\in \omega}I_m$ by $c_1<c_2<\cdots$. Now $\frac{f(|(C\setminus \cup_{m\in \omega}I_m)\cap [0,c_m-1]|)}{f(g(c_m))}=\frac{f(|(C\setminus \cup_{m\in \omega}I_m)\cap [0,c_m-1]|)}{f(h(c_m))}\to \infty$ as $m\to \infty$. By Lemma \ref{LS1}, $C\setminus \cup_{m\in \omega}I_m\in \zgf$. Next we prove that $\cup_{m\in \omega}I_m\in \zgf$.  For $k\in \omega$ choose a $m_0>\geq 2$ satisfying $k_{m_0}\leq k<k_{m_0+1}$. For all $m\geq 2$ as $\frac{f(g(k_m+\lfloor \frac{f(g(k_m))}{2^m}\rfloor))}{f(g(k_m))}>m$, $\lfloor \frac{\lfloor f(g(k_{m_0}))}{2^{m_0}} \rfloor>1$. Again, $f(g(k_m+\lfloor \frac{f(g(k_m))}{2^m}\rfloor))>mf(g(k_m))\geq 2f(g(k_m))$, i.e., $f(g(k_{m+1}))>2f(g(k_m))$, i.e., $\frac{f(g(k_{m+1}))}{2^{m+1}}>\frac{f(g(k_m))}{2^m}$ and consequently $|I_{m+1}|\geq |I_m|$. We now have

\[\frac{f(|\cup_{m\in \omega}I_m\cap [0,k-1]|)}{f(g(k))}\leq \frac{f(|\Sigma_{m\leq m_0}I_m|)}{f(g(k_{m_0}))}\leq \frac{f(m_0|I_{m_0}|)}{f(g(k_{m_0}))}\leq \frac{|I_{m_0}|f(m_0)}{f(g(k_{m_0}))}\] \[\leq (\lfloor\frac{f(g(k_{m_0}))}{2^{m_0}}\rfloor +1)\frac{f(m_0)}{f(g(k_{m_0}))}<2\frac{f(g(k_{m_0}))}{2^{m_0}}\frac{f(m_0)}{f(g(k_{m_0}))}=\frac{f(m_0)}{2^{m_0-1}}\leq \frac{m_0f(1)}{2^{m_0-1}}\to 0.\]

Thus $\lim_{k\to \infty}\frac{f(|\cup_{m\in \omega}I_m\cap [0,k-1]|)}{f(g(k))}=0$ and so $\cup_{m\in \omega} I_m\in \zgf$. Consequently, $C= (C\setminus \cup_{m\in \omega} I_m)\cup (\cup_{m\in \omega} I_m)\in \zgf$. Hence, $\zgf=\zhf$ and $h\in \Sfg$.

\ep

\begin{Cor}
\label{CS1}
Let $g\in H^{\uparrow}$ and $f$ be an unbounded modulus function. If there are a $l\in \omega$ and an increasing sequence $\{k_m\}_{m\in \omega}$ satisfying $\lim_{m\to \infty}\frac{f(g(k_m+l))}{f(g(k_m))}=\infty$, then $|\Sfg|=\ck$.
\end{Cor}
\bp
We show that for every $M>0$ and $\varepsilon>0$, $\frac{f(g(k+\lfloor \varepsilon f(g(k))\rfloor))}{f(g(k))}>M$ for infinitely many $k\in \omega$. Fix a $M>0$ and a $\varepsilon>0$. Choose a $m_0\in \omega$ satisfying $\varepsilon f(g(k_m))>l$ and $\frac{f(g(k_m+l))}{f(g(k_m))}>M$ for each $m>m_0$. We now have
\[\frac{f(g(k_m+\lfloor \varepsilon f(g(k_m))\rfloor))}{f(g(k_m))}\geq \frac{f(g(k_m+\lfloor l\rfloor))}{f(g(k_m))}=\frac{f(g(k_m+l))}{f(g(k_m))}>M\]
By Theorem \ref{TS1}, $|\Sfg|>1$. Again from Proposition \ref{PS1}, $|\Sfg|=\ck$.
\ep

We now present an example of a $g\in H^\uparrow$ and a modulus function $f$ for which $|\Sfg|=\ck$.
\begin{Ex}
\label{ES1}
Suppose that $g:\omega\to \omega$ is a function given by $g(0)=0, g(1)=1$ and $g(k)=2^{(m+1)!}$ for $k\in [2^{m!},2^{(m+1)!})$. Take the modulus function $f(x)=log(1+x)$. Clearly, $g(k)\to \infty$. Let $k_m=2^{m!}-1$ for $m\in \omega$. Now $\frac{k_m}{g(k_m)}=\frac{2^{m!}-1}{2^{m!}}=1-\frac{1}{2^{m!}}\to 1$. Therefore, $\frac{k}{g(k)}\not\to 0$ and $g\in H^\uparrow$. Let $\varepsilon>0$. Choose a $m_0\in \omega$ satisfying $\lfloor \varepsilon f(g(k_m))\rfloor\geq 1$ for all $m\geq m_0$. We now have
\[\frac{f(g(k_m+\lfloor \varepsilon f(g(k_m))\rfloor))}{f(g(k_m))}\geq \frac{f(g(k_m+1))}{f(g(k_m))}=\frac{f(g(2^{m!}))}{f(g(2^{m!}-1))}=\frac{log(1+2^{(m+1)!})}{log(1+2^{m!})}>\frac{log(2^{(m+1)!})}{log(2^{m!}+2^{m!})}\] \[=\frac{(m+1)!log(2)}{(m!+1)log(2)}>\frac{(m+1)}{2}\]
Therefore, $\lim_{m\to \infty}\frac{f(g(k_m+\lfloor \varepsilon f(g(k_m))\rfloor))}{f(g(k_m))}=\infty$.
By Theorem \ref{TS1}, $|\Sfg|>1$. Hence, by Proposition \ref{PS1}, $|\Sfg|=\ck$.
\end{Ex}

\section{Further results on $\zgf$ ideals}
In this section we now present further observations on $\zgf$ ideals.
Adapting the technique of \cite[Lemma 3.1]{mbpdmfjs} (see also \cite[Proposition 7]{ak-1}), we obtain the following result.

\begin{Prop}
\label{PD1}
Let $g\in H^\uparrow$ and $f$ be an unbounded modulus function. Let $k_0=0$ and $k_m=\min\{k\in \omega:f(g(k))\geq 2^m\}$, $m\geq 1$. Then
\[\zgf=\{C\subseteq \omega: \lim_{m\to \infty}\frac{f(|C\cap [k_m,k_{m+1})|)}{f(g(k_m))}=0\}\]
\end{Prop}
\bp
Suppose that $C\in \zgf$. For $\varepsilon>0$ choose a $k_0\in \omega$ for which $\frac{f(|C\cap [0,k-1]|)}{f(g(k))}<\frac{\varepsilon}{4}$ for each $k\geq k_0$. Also, choose a $k_{m_0}\in \omega$ for which $\frac{f(1)}{f(g(k_m))}<\frac{\varepsilon}{2}$ for each $k_m\geq k_{m_0}$. From the definition of $k_m$ it follows that $f(g(k_m))\geq 2^m$, $f(g(k_{m+1}-1))< 2^{m+1}$ and $\frac{f(g(k_{m+1}-1))}{f(g(k_m))}<2$. Now for each $m\in \omega$ and $k_m>\max\{k_0,k_{m_0}\}$ we obtain
\[\frac{f(|C\cap [k_m,k_{m+1})|)}{f(g(k_m))}\leq \frac{f(|C\cap [0,k_{m+1}-1]|)}{f(g(k_m))}\leq \frac{f(|C\cap [0,k_{m+1}-2]|+1)}{f(g(k_m))}\] \[\leq \frac{f(|C\cap [0,k_{m+1}-2]|)}{f(g(k_{m+1}-1))}\frac{f(g(k_{m+1}-1))}{f(g(k_m))}+\frac{f(1)}{f(g(k_m))}\] \[<\frac{\varepsilon}{4}2+\frac{\varepsilon}{2}=\varepsilon.\]

Hence, $\lim_{m\to \infty}\frac{f(|A\cap [k_m,k_{m+1}|)}{f(g(k_m))}=0$.

For the converse part let $C\subseteq \omega$ for which $\lim_{m\to \infty}\frac{f(|C\cap [k_m,k_{m+1})|)}{f(g(k_m))}=0$. Let $\varepsilon>0$ and choose a $m_0\in \omega$ satisfying $\frac{f(|C\cap [k_m,k_{m+1})|)}{f(g(k_m))}<\frac{\varepsilon}{4}$ for each $m\geq m_0$. Also, choose a $m_1\geq m_0$ so that $\frac{f(|C\cap [0,k_{m_0})|)}{f(g(k_{m_1}))}<\frac{\varepsilon}{4}$. Fix a $k> k_{m_1}$ and choose a $l(>m_1)\in \omega$ with $k\in [k_l,k_{l+1})$.
Observe that for $r\in \omega$ if $k_{l-r}<k_l$, then $\frac{f(g(k_{l-r}))}{f(g(k_l))}\leq\frac{1}{2^{r-1}}$.
\[\frac{f(|C\cap [0,k-1)|)}{f(g(k))}\leq \frac{f(|C\cap [0,k-1)|)}{f(g(k_l))}=\frac{f(|C\cap [0,k_{m_0}-1]|)}{f(g(k_l))}+\frac{f(|C\cap [k_{m_0},k_{m_0+1})|)}{f(g(k_l))}+\] \[+\cdots +\frac{f(|C\cap [k_{l-1},k_l)|)}{f(g(k_l))}+ \frac{f(|C\cap [k_l,k_{l+1})|)}{f(g(k_l))}\] \[\leq \frac{f(|C\cap [0,k_{m_0}-1]|)}{f(g(k_{m_1}))}+\frac{f(|C\cap [k_{m_0},k_{m_0+1})|)}{2^{l-m_0-1}f(g(k_{m_0}))}+\cdots+ \frac{f(|C\cap [k_{l-1},k_l)|)}{f(g(k_{l-1}))}+ \frac{f(|C\cap [k_l,k_{l+1})|)}{f(g(k_l))}\] \[\leq \frac{\varepsilon}{4}+(\frac{1}{2^{l-m_0-1}}+\cdots+1)\frac{\varepsilon}{4}+ \frac{\varepsilon}{4}<\frac{\varepsilon}{4}+2\frac{\varepsilon}{4}+ \frac{\varepsilon}{4}=\varepsilon.\]

Therefore, $\lim_{k\to \infty}\frac{f(|C\cap [0,k-1)|)}{f(g(k))}=0$. Hence, $C\in \zgf$.
\ep

In line with \cite[Theorem 4.11]{jt} (see also \cite[Theorem 3.7]{kbpdak}) we have the following.
\begin{Th}
\label{TD1}
Let $g\in H^\uparrow$ and $f$ be an unbounded modulus function. $\zc_g(f)$ is a generalized density ideal and $\zgf=\{C\subseteq \nb:\lim_{m\to \infty}\phi_m(C)=0\}$, where $\{\phi_m\}_{m\in \omega}$ is a sequence of submeasures concentrated on intervals $[k_m,k_{m+1})$ given by
\[\phi_m(C)=\frac{f(|C\cap [k_m,k_{m+1})|)}{f(g(k_m))},\]
and $\{k_m\}_{m\in \omega}$ is a sequence defined in Proposition \ref{PD1} and $C\subseteq \omega$.
\end{Th}

Note that $\limsup_{m\to\infty}\phi_m(\omega)>0$ as $\omega\not\in \zc_g(f)$.


\begin{Prop}
\label{PD2}
Let $g\in H^\uparrow$ and $f$ be an unbounded modulus function and $\zc_g(f)$ be a generalized density ideal generated by a sequence of submeasures $\{\phi_k\}_{k\in \omega}$ given in Theorem \ref{TD1}. The following statements are equivalent.\\
$(1)$ $\sup_{k\in \omega} \phi_k(\omega)<\infty$.\\
$(2)$ The sequence $\{\frac{f(|(f\circ g)^{-1}([2^k,2^{k+1}))|)}{2^k}\}_{k\in \omega}$ is bounded.
\end{Prop}
\bp
By Theorem \ref{TD1}, $\zgf=\{C\subseteq \nb:\lim_{m\to \infty}\phi_m(C)=0\}$,  where $\{\phi_m\}_{m\in \omega}$ is a sequence of submeasures concentrated on intervals $[k_m,k_{m+1})$ given by
$\phi_m(C)=\frac{f(|C\cap [k_m,k_{m+1})|)}{f(g(k_m))}$, $C\subseteq \omega$. From the monotocity of $g$ and $f$, we have $|\omega\cap [k_m,k_{m+1})|=|(f\circ g)^{-1}([2^m,2^{m+1}))|$. Now
\[\phi_m(\omega)=\frac{f(|\omega\cap [k_m,k_{m+1})|)}{f(g(k_m))}=\frac{f(|(f\circ g)^{-1}([2^m,2^{m+1}))|)}{f(g(k_m))}.\]

Clearly, $\sup_{k\in \omega} \phi_k(\omega)<\infty$ implies that $\sup_{k\in \omega}\frac{f(|(f\circ g)^{-1}([2^k,2^{k+1}))|)}{2^k}<\infty$. Hence, $(2)$ holds.
Again if $(2)$ holds, then $\sup_{k\in \omega} \phi_k(\omega)<\infty$. Hence, $(1)$ holds.
\ep

\begin{Prop}
\label{PD3}
Let $g\in H^\uparrow$ and $f$ be an unbounded modulus function and $\zc_g(f)$ be a generalized density ideal generated by a sequence of submeasures $\{\phi_k\}_{k\in \omega}$ given in Theorem \ref{TD1}. Suppose that for each $M>0$ there exists a $L>0$ satisfying $\frac{f(g(k+\lfloor Lf(g(k)) \rfloor))}{f(g(k))}>M$ for all but finitely many $k\in \omega$. Then $\sup_{k\in \omega} \phi_k(\omega)<\infty$.
\end{Prop}
\bp
Suppose that $M=2$. By the given condition, there are $L>0$ and $k_0\in \omega$ for which $\frac{f(g(k+\lfloor Lf(g(k)) \rfloor))}{f(g(k))}>2$ for all $k\geq k_0$. Assume that $\alpha=2L$. We prove that
$f(|(f\circ g)^{-1}([2^k,2^{k+1}))|)\leq 2^kf(\alpha)$ for all $k\in \omega$ satisfying $2^k\geq f(g(k_0))$. Let $k\in \omega$ and let $m=\min(f\circ g)^{-1}([2^k,2^{k+1}))$. For $(f\circ g)^{-1}([2^k,2^{k+1}))=\emptyset$, the case is trivial. Clearly, $m\geq k_0$. Therefore,
\[2<\frac{f(g(m+\lfloor Lf(g(m)) \rfloor))}{f(g(m))}\leq \frac{f(g(m+\lfloor L2^{k+1} \rfloor))}{2^k}\leq \frac{f(g(m+\lfloor 2^k\alpha \rfloor))}{2^k}.\]
We obtain that $f(g(m+\lfloor 2^k\alpha \rfloor))\geq 2^{k+1}$. Consequently, $|(f\circ g)^{-1}([2^k,2^{k+1}))|\leq  \lfloor 2^k\alpha \rfloor\leq 2^k\alpha$ which implies $f(|(f\circ g)^{-1}([2^k,2^{k+1}))|)\leq f(2^k\alpha)\leq 2^kf(\alpha)$. Thus $\frac{f(|(f\circ g)^{-1}([2^k,2^{k+1}))|)}{2^k}\leq f(\alpha)$. By Proposition \ref{PD2}, $\sup_{k\in \omega} \phi_k(\omega)<\infty$.
\ep

Adapting the technique of \cite[Proposition 1]{ak-2} we have the following observation.

\begin{Prop}
\label{PD4}
Let $g\in H^\uparrow$ and $f$ be an unbounded modulus function and $\zc_g(f)$ be a generalized density ideal generated by a sequence of submeasures $\{\phi_k\}_{k\in \omega}$ given in Theorem \ref{TD1}. The following statements are equivalent.\\
$(1)$ $\sup_{k\in \omega} \phi_k(\omega)<\infty$.\\
$(2)$ The sequence $\{\frac{f(k)}{f(g(k))}\}_{k\in \omega}$ is bounded.
\end{Prop}
\bp
$(1)\Rightarrow (2)$. By \cite[Theorem 4.11]{jt}, $\zc_g(f)$ is generated by a sequence $\{\phi_m\}_{m\in \omega}$ of submeasures, where $\phi_m(C)=\frac{f(|C\cap [k_m,k_{m+1})|)}{f(g(k_m))}$ for some $C\subseteq \omega$ and $k_{m+1}=\min\{k\in \omega:f(g(k))\geq 2f(g(k_m))\}$. Choose a $M>0$ such that $\phi_m(\omega)<M$ for each $m\in \omega$. Now for $k\in \omega$, $k_m<k\leq k_{m+1}$ for some $m\in \omega$ and in the view of the fact that $f(k_{m+1})\leq f(k_{m+1}-k_m)+f(k_m)$ one obtains
\[\frac{f(k)}{f(g(k))}\leq \frac{f(k_{m+1})}{f(g(k_m))}=\frac{f(k_0)}{f(g(k_m))}+\Sigma^m_{j=0}\frac{f(k_{j+1})-f(k_j)}{f(g(k_m))}\leq \frac{f(k_0)}{f(g(k_0))}+\Sigma^m_{j=0}\frac{f(k_{j+1})-f(k_j)}{2^{m-j}f(g(k_j))}\] \[\leq \frac{f(k_0)}{f(g(k_0))}+\Sigma^m_{j=0}\frac{f(k_{j+1}-k_j)}{2^{m-j}f(g(k_j))}= \frac{f(k_0)}{f(g(k_0))}+\Sigma^m_{j=0}\frac{\phi_j(\omega)}{2^{m-j}}\] \[\leq \frac{f(k_0)}{f(g(k_0))}+2M.\]

This shows that $\{\frac{f(k)}{f(g(k))}\}_{k\in \omega}$ is bounded.

$(2)\Rightarrow (1)$. Choose a $M>0$ such that $\frac{f(k)}{f(g(k))}< M$ for each $k\in \omega$. Clearly, $f(g(0))>0$ and also $f(g(k_{m+1}-1))<2f(g(k_m))$. It now follows that
\[\phi_m(\omega)=\frac{f(|\omega\cap [k_m,k_{m+1})|)}{f(g(k_m))}=\frac{f(k_{m+1}-k_m)}{f(g(k_m))}\leq\frac{f(k_{m+1})}{f(g(k_m))}< \frac{2f(k_{m+1})}{f(g(k_{m+1}-1))}\] \[\leq \frac{2f(k_{m+1}-1)+2f(1)}{f(g(k_{m+1}-1))}<2M+\frac{2f(1)}{f(g(0))}.\]
Therefore, $\sup_{m\in \omega}\phi_m(\omega)<\infty$.
\ep

\begin{Prop}
\label{PD4-2}
Let $g\in H^\uparrow$ and $f$ be an unbounded modulus function and $\zc_g(f)$ be a generalized density ideal generated by a sequence of submeasures $\{\phi_k\}_{k\in \omega}$ given in Theorem \ref{TD1}. If the sequence $\{\frac{k}{g(k)}\}_{k\in \omega}$ is bounded, then $\sup_{k\in \omega} \phi_k(\omega)<\infty$.
\end{Prop}
\bp
Choose an $M\in \omega$ for which $\frac{k}{g(k)}\leq M$ for each $k\in \omega$, i.e., $k\leq Mg(k)$. As $f$ is a modulus function, $f(k)\leq f(Mg(k))\leq Mf(g(k))$. Thus for each $k\in \omega$ $\frac{f(k)}{f(g(k))}\leq M$. Hence, by Proposition \ref{PD4}, $\sup_{k\in \omega} \phi_k(\omega)<\infty$.
\ep

\begin{Prop}
\label{PD4-1}
Let $g\in H^\uparrow$ and $f$ be an unbounded modulus function. If the sequence $\{\frac{f(k)}{f(g(k))}\}_{k\in \omega}$ is bounded, then $\zf\subseteq \zgf$.
\end{Prop}
\bp
Choose an $M>0$ satisfying $\frac{f(k)}{f(g(k))}\leq M$ for each $k\in \omega$. Let $C\in \zf$. For $k\in \omega$ we obtain
\[\frac{f(|C\cap [0,k-1]|)}{f(g(k))}=\frac{f(|C\cap [0,k-1]|)}{f(k)}\frac{f(k)}{f(g(k))}\leq M\frac{f(|C\cap [0,k-1]|)}{f(k)}\]
As $C\in \zf$, $\limsup_{k\to \infty}\frac{f(|C\cap [0,k-1]|)}{f(k)}=0$. Therefore, $\lim_{k\to \infty}\frac{f(|C\cap [0,k-1]|)}{f(g(k))}=0$ and $C\in \zgf$. Hence, $\zf\subseteq \zgf$.
\ep

The following example shows that the converse of Proposition \ref{PD4-2} is not true.
\begin{Ex}
\label{EEU3}
Suppose that $g:\omega\to \omega$ is a function defined by $g(k)=1$ for $k\leq 4$ and $g(k)=2^m$ for $4^m<k\leq 4^{m+1}$, $m>1$ and $f(x)=log(1+x)$. Choose $k_m=4^{m+1}$, $m\in \omega$. Then $\frac{k_m}{g(k_m)}=\frac{4^{m+1}}{2^m}\to \infty$ as $m\to \infty$. Now for $k\in \omega$ choose a $m\in \omega$ with $4^m<k\leq 4^{m+1}$ and $g(k)=2^m$.
\[\frac{f(k)}{f(g(k))}=\frac{log(1+k)}{log(1+2^m)}\leq \frac{log(1+4^{m+1})}{log(1+2^m)}<\frac{2log(1+2^{m+1})}{log(1+2^m)}\leq 4\]
Therefore, the sequence $\{\frac{f(k)}{f(g(k))}\}_{k\in \omega}$ is bounded. By Proposition \ref{PD4}, $\sup_{k\in \omega} \phi_k(\omega)<\infty$.
\end{Ex}

We now present an example which shows that the converse of Proposition \ref{PD3} is not true.
\begin{Ex}
\label{EEU}
Suppose that $g:\omega\to \omega$ is a function defined by $g(0)=0$ and $g(k)=(m+1)!$ for all $k\in [m!,(m+1)!)$, $k,m\in \omega$. Let $f(x)=log(1+x)$. For $k\in \omega$ choose a $m\in \omega$ with $k\in [m!,(m+1)!)$. Now $\frac{f(k)}{f(g(k))}=\frac{f(k)}{f((m+1)!)}< 1$. By Proposition \ref{PD4}, $\sup_{k\in \omega}\phi_k(\omega)<\infty$. Let $M=2$ and $L>0$. Define $k_m=m!$ for $m\in \omega$. Choose a $m_0\in \omega$ with $m_0>L$. For $m>m_0$ we have $Lf(g(k_m))=Llog(1+(m+1)!)<L(m+1)!$ and $k_m+\lfloor Lf(g(k_m))\rfloor< (m+1)!+L(m+1)!=(L+1)(m+1)!$, i.e., $k_m+\lfloor Lf(g(k_m))\rfloor<(m+2)!$.
\[\frac{f(g(k_m+\lfloor Lf(g(k_m))\rfloor))}{f(g(k_m))}\leq \frac{f(g((m+2)!-1))}{f(g(k_m))}\leq \frac{f(g((m+2)!-1))}{f(g(m!))}=\frac{log(1+(m+2)!)}{log(1+(m+1)!)}\leq 2.\]

Hence, $\frac{f(g(k+\lfloor Lf(g(k))\rfloor))}{f(g(k))}\leq 2$ for infinitely many $k\in \omega$.
\end{Ex}


%

Note that $\zc_{id_\omega}(f)$ is not a simple density ideal (see \cite[Example 2.3]{kbpdak}). We now present an example of an ideal $\zgf$ which is not a simple density ideal for a nonidentity function $g$.

\begin{Ex}
\label{EEC}
Suppose that $g:\omega\to \omega$ is a function defined by $g(k)=1$ for $k\leq 4$ and $g(k)=2^m$ for $4^m<k\leq 4^{m+1}$, $m>1$ and let $f(x)=log(1+x)$. Let $h\in H$. We show that $\zgf\neq \zh$. We now have two possible cases.

Case 1: There exists $\alpha>0$ and $m_0\in \omega$ for which $1+h(k)>k^\alpha$ for all $k\geq m_0$.

Let $\alpha<1$ and $m_0>1$.  Consider the set $C\subseteq \omega$, where $|C\cap [0,k-1]|=\lfloor k^\frac{\alpha}{2}\rfloor$. Clearly, $C\in \zh$. For $k\in \omega
$ $4^m<k\leq 4^{m+1}$ and $g(k)=2^m$ for some $m\in \omega$, we now have
\[\frac{f(|C\cap [0,k-1]|)}{f(g(k))}=\frac{log(1+|C\cap [0,k-1]|)}{log(1+g(k))}\geq \frac{log(k^\frac{\alpha}{2})}{log(1+g(k))}=\frac{\frac{\alpha}{2}log(k)}{log(1+2^m)}\] \[>\frac{\frac{\alpha}{2}log4^m}{log2^{m+1}}=\frac{\alpha m}{m+1}.\]

Therefore, $\lim_{k\to \infty}\frac{f(|C\cap [0,k-1]|)}{f(g(k))}\geq \alpha$ ($m\to \infty$ as $k\to \infty$). Hence, $C\not \in \zgf$ and $\zgf\neq \zh$.

Case 2: For every $\alpha>0$ and $m\in \omega$ there exists a $k>m$ for which $1+h(k)\leq k^\alpha$.

Choose an increasing sequence $\{k_m\}$ from $\omega$ satisfying $1+h(k_m)\leq k^\frac{1}{m+1}_m$ and $k^\frac{1}{m+1}_m<k^\frac{1}{m+2}_{m+1}$ for every $m\in \omega$. Consider the set $C\subseteq \omega$, where $|C\cap [0,k_m-1]|=\lfloor k^\frac{1}{m+1}_m-1\rfloor$ and $|C\cap [0,k-1]|\leq k^\frac{1}{m+1}-1$ for each $k$, with $k_m\leq k<k_{m+1}$. Clearly, $C\not\in \zh$. For $k\in \omega$ choose $m,l\in \omega$ with $k_m\leq k<k_{m+1}$ and $4^l<k\leq 4^{l+1}$. Then
\[\frac{f(|C\cap [0,k-1]|)}{f(g(k))}\leq \frac{log(k^\frac{1}{m+1})}{log(1+g(k))}=\frac{1}{m+1}\frac{log(k)}{log(1+g(k))}\leq \frac{1}{m+1}\frac{log(4^{l+1})}{log(1+2^l)}\] \[<\frac{1}{m+1}\frac{l+1}{l}\frac{log(4)}{log(2)}<\frac{2}{m+1}(1+\frac{1}{l})<\frac{4}{m+1}.\]
As $k\to \infty$ implies that $m\to \infty$, $\lim_{k\to \infty}\frac{f(|C\cap [0,k-1]|)}{f(g(k))}=0$. Hence, $C\in \zgf$ and $\zgf\neq \zh$.
\end{Ex}


%
%

\subsection{Some observations on $\zgf$ related to increasing-invariance}
In this section we investigate $\zgf$ ideals in the context of the notion of increasing invariance which was introduced in \cite{ak-2}. Let $g\in G$ and $f$ be an unbounded modulus function. Then $\zgf$ is increasing invariant \cite[Page 13]{jt}.
%
Following the same argument as in \cite[Proposition 5]{ak-2} we have the next result.
\begin{Prop}
\label{PD6}
There exists an $\sqsubseteq$-antichain of size $\ck$ among $\EU$ ideals that are not $\zgf$ for any $g\in G$ and any modulus function $f$.
\end{Prop}
\bp
Suppose that $\fc$ is a collection of pairwise almost disjoint infinite subsets of $\omega$ of cardinality $\ck$. Let $k_0=1$ and $k_{j+1}=k_j+j!+(j+1)!$, $j\in \omega$. Consider a sequence $\{\mu_j\}_{j\in \omega}$ of probablity measures having a finite support $D_j=\{k_j,k_j+1,\cdots, k_j+j!-1\}$ satisfying $\mu_j(\{i\})=\frac{1}{j!}$ for $i\in D_j$ and $\mu_j(\{i\})=0$ otherwise. For each $L\in \fc$, let $\ic_L=\Exh(\sup_{l\in L} \mu_l)$. $\ic_L$ is an $\EU$ ideal. Now $\ic_L$ is not a $\zgf$ ideal for any $g\in G$ and for any modulus function $f$ because $\zgf$ is increasing invariant but $\ic_L$ is not increasing invariant. Also for $L,K\in \fc$, $\ic_L\not\sqsubseteq \ic_K$ (see \cite[Proposition 5]{ak-2}).
\ep

We know that for a $\beta\in (0,1)$, $f(x)=x^\beta$, $x\in \rb^{+}\cup \{0\}$ and for any $g\in G$ $\zgf=\zg$. Thus from \cite[Proposition 6]{ak-2} it follows that
\begin{Prop}
\label{PD7}
There exists an $\sqsubseteq$-antichain of size $\ck$ among $\EU$ $\zgf$ ideals for some $g\in G$ and for some modulus function $f$.
\end{Prop}

From \cite[Theorem 6]{ak-1} and the fact that $\zgf$ is increasing invariant we obtain.
\begin{Prop}
\label{PD10}
Let $g\in G$ and $f$ be an unbounded modulus function. Let $A\in H(\zgf)$ and $A=\{a_0<a_1<\cdots\}$. The function $h:\omega\to A$ given by $h(k)=a_k$ is a bijection for which $\zgf|A\cong \zgf$ holds.
\end{Prop}

We now present some examples of ideals that are not $\zgf$ for any $g\in G$ and for any modulus function $f$.
\begin{Ex}
\label{EEU4}
Consider the ideal $\ic=\Exh(\phi)$ defined in \cite[Example 3.2]{mbpdmfjs}, where $\phi=\sup_{m\in \omega}\mu_m$ and $\mu_m$ is a measure with the support $\omega\cap [2^m,2^m+m)$ satisfying $\mu_m(\{k\})=\frac{1}{m}$ for all $k\in [2^m,2^m+m)$. Let $C=\cup_{m\geq 1}[2^m,2^m+m)\cap \omega$ and $D=\cup_{m\geq 1}[2^m-m,2^m)\cap \omega$. Then $C\not\in \ic$ and $D\in \ic$. Note that for any $k\in \omega$ $|C\cap [0,k-1]|\leq |D\cap [0,k-1]|$. Therefore $\ic$ is not increasing invariant. Since $\zgf$ is increasing invariant, $\ic$ is not $\zgf$ for any $g\in G$ and for any modulus function $f$.

\end{Ex}


\begin{Ex}
\label{EEU5}
Consider the ideal $\ic=\Exh(\phi)$ defined in \cite[Example 3.3]{mbpdmfjs}, where $\phi=\sup_{m\in \omega}\mu_m$ and $\mu_m$ is a measure with a support $\omega\cap [2^m,2^m+m^2)$ satisfying $\mu_m(\{k\})=\frac{1}{m}$ for all $k\in [2^m,2^m+m^2)$. Then $\ic$ is not increasing invariant and hence it is not $\zgf$ for any $g\in G$ and for any modulus function $f$.

\end{Ex}

Finally, we present the following observations on $\zf$.
\begin{Lemma}
\label{LZ}
Let $f$ be an unbounded modulus function. Then $\zf\subseteq \zc$.
\end{Lemma}
\bp
Suppose that $C\in \zf$ and $\varepsilon>0$. Choose a $M\in \omega$ with $\frac{1}{M}<\varepsilon$. We have $\limsup_{k\to \infty} \frac{f(|C\cap [0, k-1]|)}{f(k)}=0$. Choose a $k_0\in \omega$ satisfying $\frac{f(|C\cap [0, k-1]|)}{f(k)}<\frac{1}{M}$ for all $k\geq k_0$. Now $f(M|C\cap [0,k-1]|)\leq Mf(|C\cap [0,k-1]|)<f(k)$. Consequently, for all $k\geq k_0$ $M|C\cap [0,k-1]|<k$ (as $f$ is increasing), i.e., $\frac{|C\cap [0,k-1]|}{k}<\frac{1}{M}<\varepsilon$. Therefore, $C\in \zc$. Hence, $\zf\subseteq \zc$.
\ep

If $\ic$ is an increasing-invariant $\EU$ ideal, then $\zc\subseteq \ic$ \cite[Lemma 1]{ak-2}. From Lemma \ref{LZ} and \cite[Lemma 1]{ak-2} it follows that
\begin{Lemma}
\label{LZ-2}
Let $f$ be an unbounded modulus function. If $\ic$ is an increasing-invariant $\EU$ ideal, then $\zf\subseteq \ic$.
\end{Lemma}

\noindent{\textbf{Acknowledgement:}} The first author as PI is thankful to SERB(DST) for the CRG project (Project No. CRG/2022/000264)  during the tenure of which this work has been done. The second author (as fellow) is thankful to NBHM for the research project (Project No. 02011/9/2022/NBHM(R.P)/R$\&$D II/10378) during the tenure of which this work has been done.\\

\noindent{\textbf{Author Contributions:}} Both authors wrote the paper and reviewed it.\\

\noindent{\textbf{Data Availability:}} Not applicable.\\

\noindent{\textbf{Declarations}}\\
\noindent{\textbf{Conflict of Interest:}} The authors state that there is no conflict of interest.

\end{document}